\documentclass[11pt]{article}
\usepackage{amssymb}
\usepackage{amsmath}
\usepackage{amsbsy}
\usepackage{bm}
\newtheorem{lemma}{Lemma}[section]
\newtheorem{proposition}{Proposition}[section]
\newtheorem{theorem}{Theorem}[section]

\def\proclaim#1{\par \bigskip\noindent {\bf #1}\bgroup\it\ }
\def\endproclaim{\egroup\par\bigskip}
\def\proof{\par\noindent{\bf Proof.} \;}

\newbox\TempBox \newbox\TempBoxA
\newcommand{\non}{\nonumber }
\def\pr{\textsf{P}} % the symbol P for probability used the sans serif letter
\def\ep{\textsf{E}} % the symbol E for expectation used the sans serif letter
\def\Cov{\textsf{Cov}} % the symbol Cov for covariance used the sans serif letter
\def\Var{\textsf{Var}} % the symbol Var for covariance used the sans serif letter

\def\bk#1{\bm #1}
\def\text#1{\mbox{\rm #1}}

\def\underwiggle 1{
\ifmmode\setbox\TempBox=\hbox{$ 1$}\else\setbox\TempBox=\hbox{
1}\fi \setbox\TempBoxA=\hbox to \wd\TempBox{\hss\char'176\hss}
\rlap{\copy\TempBox}\smash{\lower9pt\hbox{\copy\TempBoxA}} }

\begin{document}

%\renewcommand{\theequation}
%{\arabic{section}.\arabic{equation}}
%\baselineskip=22pt

%$ $ \vskip 0.8in

%\Large

\title{\huge \bf Maximal inequalities and a law of the iterated logarithm for negatively associated
 random fields\footnote{ Research supported by National Natural Science Foundation of
China  } }
\author  {{\sc  Li-Xin ZHANG}
\\
{\em Department of Mathematics, Zhejiang University, Hangzhou
310027, China} }

%\date{\today}
\date{ }

\maketitle
\begin{abstract} The exponential inequality of the maximum partial
sums is a key to establish the law of the iterated logarithm of
negatively associated random variables. In the one-indexed random
sequence case, such inequalities for negatively associated random variables
are established by Shao (2000) by
using his comparison theorem between negatively associated and
independent random variables. In the multi-indexed random field
case, the comparison theorem fails. The purpose of this paper is
to establish the Kolmogorov exponential inequality as well a moment inequality
 of the maximum
partial sums of a negatively associated random field via a
different method. By using these inequalities, the sufficient and
necessary condition for the law of the iterated logarithm of a
negatively associated random field to hold is obtained.


\vskip2mm {\bf Abbreviated Title:} LIL for NA Fields

\vskip2mm


\end{abstract}

%\end{document}

\section{  Introduction and the law of iterated logarithm}
\setcounter{equation}{0}

 A finite family of random variables $\{X_i; 1\le i\le n\}$ is
said to be negatively associated if for every pair of disjoint
subsets $A$ and $B$ of $\{1,2,\cdots,n\}$,
\begin{equation}\label{eq1.1} \Cov\{f(X_i; i\in A), g(X_j; j\in B)\}\le
0 \end{equation}
 whenever $f$ and $g$ are coordinatewise
non-decreasing and the covariance exists. An infinite family is
negatively associated if every finite subfamily is negatively
associated. The concept of the negative association was introduced
by Alam and Saxena (1981) and Joag-Dev and Proschan (1983). As
pointed out and proved by Joag-Dev and Proschan (1983), a number
of well-known multivariate distributions possess the NA property.
 Many
investigators  discuss the properties and limit theorems of negatively
 associated random variables. We refer to Joag-Dev and
Proschan (1983) for fundamental properties, Newman (1984) for the
central limit theorem, Matula (1992) for the three series theorem,
Su, et al. (1997) for the moment inequality, Roussas (1996) for
the Hoeffding inequality, and Shao (2000) for the Rosenthal-type
maximal inequality and
 the Kolmogorov exponential inequality.
 Shao and Su (1999) established the law of the iterated logarithm
for negatively associated random variables with finite variances.

\proclaim{ Theorem A}  Let $\{X_i; i\ge 1\}$ be a strictly
stationary negatively associated sequence with $\ep X_1=0$, $\ep
X_1^2<\infty$ and $\sigma^2:=\ep X_1^2+2\sum_{i=2}^{\infty}\ep
(X_1X_i)>0$. Let $S_n=\sum_{i=1}^n X_i$. Then
\begin{equation}\label{eq1.2}
\limsup_{n\to \infty} \frac{S_n}{(2n\log\log n)^{1/2} } =\sigma
\quad a.s. \end{equation}
 Here and in the squeal of this paper,
$\log x =\ln(x\vee e)$
\endproclaim

Let $\{X_{\bk n}; \bk n\in \Bbb N^d\}$ be a field of random
variables, where $d\ge 2$ is a positive integer, $\Bbb N^d$
denotes the d-dimensional lattice of  positive integers. Through
this paper, for $\bk n=(n_1,\cdots, n_d)\in \Bbb N^d$, $\bk
k=(k_1,\cdots, k_d)\in \Bbb N^d$ and $m\in \Bbb N$, we denote
$|\bk n|= n_1\cdots n_d$, $\|\bk n\|= n_1+ \cdots +n_d$, $\bk k
\bk n=(k_1 n_1,\cdots, k_d n_d)$ and $\bk k m=(k_1 m,\cdots, k_d
m)$. Also, $\bk k \le \bk n$ (resp. $ \bk k\ge \bk n$) means
$k_i\le n_i$ (resp. $k_i\ge n_i$), $i=1,2,\cdots, d$. Denote by
$S_{\bk n}=\sum_{\bk k \le \bk n} X_{\bk k}$  and $\bk
1=(1,\cdots,1)\in \Bbb N^d$. It is known that, if $\{X_{\bk n}\}$
is a field of i.i.d.r.v.s, then
\begin{equation}\label{eq1.3}
\limsup_{ \bk n \to \infty} \frac{ |S_{\bk n}| }{ (2d|\bk n|
\log\log |\bk n|)^{1/2} } =(\ep X_{\bk 1}^2)^{1/2} \quad a.s.
\end{equation}
if and only if $\ep X_{\bk 1}=0$ and $\ep X_{\bk 1}^2 \log^{d-1}
(|X_{\bk 1}|) / \log\log (|X_{\bk 1}|)<\infty$, where $ \bk n\to
\infty$ means $n_1\to \infty, \cdots, n_d\to \infty$. When $\{
X_{\bk n}; \bk n \in \Bbb N^d\}$ is a negatively associated field
of random variables, Zhang and Wen (2001a) and Zhang
and Wang (1999)
 established   the weak convergence,
the law of large numbers and the complete convergence similar to
those for fields of independent random variables. This paper is to
establish the law of the iterated logarithm similar to (1.3)
 for a negatively associated random field.

\begin{theorem}\label{th2} Let $d\ge 2$ be a positive integer, and $\{X_{\bk
n}; \bk n \in \Bbb N^d\}$ be a weakly stationary negatively
associated field of identically distributed random variables satisfying
\begin{equation}\label{eq1.4}
 \ep X_{\bk 1}=0 \quad \text{ and } \quad \ep X_{\bk
1}^2 \log^{d-1}( |X_{\bk 1}|) / \log\log (|X_{\bk 1}|)<\infty.
\end{equation}
 Denote by $ \Upsilon(\bk j-\bk i)= \Cov( X_{\bk j},
X_{\bk i})$ and $\sigma^2= \sum_{ \bk j \in \Bbb Z^d}\Upsilon(\bk
j)$. Then
\begin{equation}\label{eq1.5}
\limsup_{\bk n \to \infty} \frac{ S_{\bk n} }{ (2d|\bk n| \log\log
|\bk n|)^{1/2} } = \sigma \quad a.s.
\end{equation}
\end{theorem}

\begin{theorem}\label{th1} Let $d\ge 2$ be a positive integer, and $\{X_{\bk
n}; \bk n \in \Bbb N^d\}$ be a   negatively
associated field of identically distributed random variables satisfying (\ref{eq1.4}).
 Then
$$
\limsup_{\bk n \to \infty} \frac{ |S_{\bk n}| }{ (2d|\bk n| \log\log
|\bk n|)^{1/2} } \le (\ep X_{\bk 1}^2)^{1/2}  \quad a.s.
$$
\end{theorem}

The following theorem tells us that the condition (\ref{eq1.4}) is
necessary for the law of the iterated logarithm to hold.
\begin{theorem}\label{th3} Let $d\ge 1$ be a positive integer, and
$\{X_{\bk n}; \bk n \in \Bbb N^d\}$ be a negatively associated
field of identically distributed random variables. If
\begin{equation}\label{eq1.6}
\pr\left( \limsup_{\bk n\to \infty} \frac{ |S_{\bk n}|}{ (2 d |\bk
n| \log\log |\bk n|)^{1/2} }< \infty\right)>0,
\end{equation}
 then (\ref{eq1.4}) holds.
\end{theorem}

In showing  the law of the iterated logarithm, a main step is to
establish the exponential inequalities.  In the case of $d=1$,
such exponential inequalities for negatively associated random
variables are established by Shao (2000) by using a comparison
theorem between negatively associated random variables  and independent random
variables. However, if $d\ge 2$, such a comparison theorem fails for
the maximum partial sums (cf., Bulinski and Suquet 2001). In
section 2, we establish a Kolmogorov type exponential inequality
as well as a moment inequality
of the maximum partial sums of a negatively associated field via a
different method. Theorems \ref{th2}-\ref{th3} are proved in
Section 3.

%%%%%%%%%%%%%%%%%%%%%%%%%%%%%%%%%%
\section{Moment inequalities and exponential inequalities}
\setcounter{equation}{0} First, we have the following moment
inequalities and exponential inequalities for the partial sums.

\begin{lemma}\label{lem2.1}   Let $p\ge 2$ and let
$\{ Y_{\bk k}; \bk k \le \bk n \}$ be a negatively associated
field of random variables with $\ep Y_{\bk k} =0$ and $\ep |Y_{\bk
k}|^p<\infty$. Then
$$ \ep  \big| \sum_{\bk i \le \bk n} Y_{\bk i} \big|^p
\le 2(15p/\ln p)^p \Big\{ \Big(\sum_{\bk i \le \bk n} \ep Y_{\bk
i}^2\Big)^{p/2} +\sum_{\bk i \le \bk n} \ep |Y_{\bk i}|^p\Big\}.
$$
\end{lemma}

\begin{lemma}\label{lem2.2}
Let $\{ Y_{\bk k}; \bk k \le \bk n \}$ be a negatively associated
field of random variables with zero means and finite second
moments. Let $T_{\bk k}= \sum_{\bk i \le \bk k} Y_{\bk i}$ and
$B_{\bk n}^2 = \sum_{\bk k \le \bk n} \ep Y_{\bk k}^2$. Then for
all $x>0$ and $a>0$,
\begin{equation}\label{eq2.1}
 \pr (T_{\bk
n} \ge x) \le \pr \big( \max_{ \bk k \le \bk n } Y_{\bk k} > a
\big) +\exp \Big(-\frac{x^2}{2(ax+B_{\bk n}^2)}\Big).
\end{equation}
\end{lemma}

\begin{lemma}\label{lem2.3} Let $\{ Y_{\bk k}; \bk k \le \bk n \}$ be a
negatively associated field of random variables with $\ep Y_{\bk
k}=0$ and $\ep |Y_{\bk k}|^3<\infty$. Denote by $T_{\bk k}=
\sum_{\bk i \le \bk k} Y_{\bk i}$ and $B_{\bk n}^2 = \sum_{\bk k
\le \bk n} \ep Y_{\bk k}^2$. Then for all $x>0$,
\begin{eqnarray}\label{eq2.2}
\pr(T_{ \bk n} \ge x B_{\bk n}) &\ge& \big(1-\Phi(x+1)\big) -6
B_{\bk n}^{-2} \sum_{\bk 1\le \bk i \ne \bk j \le \bk n}
|\ep(Y_{\bk i}Y_{\bk j})|
\non\\
&&-12 B_{\bk n}^{-3} \sum_{ \bk k \le \bk n} \ep |Y_{\bk k}|^3
\end{eqnarray}
where $\Phi(x)$ is the distribution of a standard normal variable.
\end{lemma}

\noindent{\bf Proofs of Lemmas \ref{lem2.1}- \ref{lem2.3}:} In the
case of $d=1$, Lemma  \ref{lem2.1} is proved by Shao (2000), and
Lemma \ref{lem2.3} is proved by Shao and Su (1999). Also,
(\ref{eq2.1}) follows from the following inequality easily:
$$\pr (T_{\bk n} \ge x) \le \pr \big( \max_{ \bk k \le \bk n } Y_{\bk k} > a \big)
+\exp \left\{\frac xa-\Big(\frac xa+\frac{B_n^2}{a^2}\Big)\ln\Big(\frac
{xa}{B_n^2}+1\Big)\right\}. $$ The later is proved by Su, et
al.(1997). Since Lemmas \ref{lem2.1}-\ref{lem2.3} do not involve
the partial order of the index set, so them are valid for $d\ge 2$
also. In fact, when $d\ge 2$, there is a one-one map $\pi: \{ \bk
k: \bk k \le \bk n\} \to \{1,2, \cdots, |n|\}$. By noting that
$\{Y_{\pi^{-1}(i)}:i=1,\cdots, |\bk n|\}$ is a negatively
associated sequence and $\sum_{i=1}^{|\bk n|}(\cdot)= \sum_{\bk k
\le \bk n}(\cdot)$, $ \max_{\bk k\le \bk n} Y_{\bk k}=\max_{i\le
|\bk n|}Y_{\pi^{-1}(i)}$,  the results follow.

In a same may, one can extend (1.2) of shao (2000) to the case of
$d\ge 2$.
\begin{lemma}\label{lem2.4} Let $\{ Y_{\bk k}; \bk k \le \bk n \}$ be
a negatively associated field and $\{ Y_{\bk k}^{\ast}; \bk k \le
\bk n \}$ be a field of independent random variables such that for
each $\bk k$, $Y_{\bk k}$ and $Y_{\bk k}^{\ast}$ have the same
distribution. Then
$$ \ep f(\sum_{\bk k\le \bk n} Y_{\bk k})
\le \ep f(\sum_{\bk k\le \bk n} Y_{\bk k}^{\ast})
$$
for any convex function $f$ on $\Bbb R$, whenever the expectations
exist.
\end{lemma}

It shall be mentioned that it is impossible to find a one-one map
$\pi: \{ \bk k: \bk k \le \bk n\} \to \{1,2, \cdots, |n|\}$ such
that $\sum_{i=1}^{|\bk m|}Y_{\pi^{-1}(i)}= \sum_{\bk k \le \bk
m}Y_{\bk k}$ for all $\bk m\le \bk n$. So, the inequalities for
maximum partial sums can not be extended directly.

For maximum partial sums, Zhang and Wen (2001a) established two
moment inequalities.

\begin{proposition}\label{prop2.1}
{\rm (Zhang and Wen 2001a)} Let $p\ge 2$, and let $\{ Y_{\bk k};
\bk k \le \bk n \}$ be a negatively associated field of random
variables with $\ep Y_{\bk k}=0$ and $\ep |Y_{\bk k}|^p<\infty$.
Suppose that $\{ \epsilon_{\bk k}; \bk k \le \bk n \}$ is a field
of i.i.d.r.v.s with $\pr(\epsilon_{\bk k} = \pm 1)=1/2$. Also,
suppose that $\{ \epsilon_{\bk k}; \bk k \le \bk n \}$ is
independent of $\{ Y_{\bk k}; \bk k \le \bk n \}$. Denote by
$T_{\bk k}= \sum_{\bk i \le \bk k} Y_{\bk i}$, $M_{\bk
n}=\max_{\bk k \le \bk n} |T_{\bk k}|$, $\widetilde T_{\bk k}=
\sum_{\bk i \le \bk k} Y_{\bk i}$, $\widetilde M_{\bk n}=\max_{\bk
k \le \bk n} |\widetilde T_{\bk k}|$ and $\|\bk X\|_p=(\ep
|X|^p)^{1/p}$. Then
\begin{equation} \label{eq2.3}
 \|M_{\bk n}\|_p\le 5 \|\widetilde M_{\bk
n}\|_p+\|M_{\bk n}\|_1 ,
\end{equation}
 and there exists a constant $A_{p,d}$ depending on $p$ and $d$  such that
\begin{equation}\label{eq2.4}
\ep | M_{\bk n} |^p \le A_{p,d} \Big\{ (\ep | M_{\bk n} | )^p
+\big( \sum_{\bk k \le \bk n} \ep Y_{\bk k}^2 \big)^{p/2}
+\sum_{\bk k \le \bk n} \ep |Y_{\bk k}|^p \Big\}.
\end{equation}
\end{proposition}

\begin{proposition}\label{prop2.2}
 {\rm (Zhang and Wen 2001a)}  Let $\{ Y_{\bk k};
\bk k \le \bk n \}$ be a strictly stationary negatively associated
field of  random variables with $\ep Y_{\bk 1}=0$ and $0<\ep
Y_{\bk 1}^2<\infty$. Denote by $T_{\bk k}= \sum_{\bk i \le \bk k}
Y_{\bk i}$. Then there exists a constant $K$, depending only on
$d$, such that
\begin{equation}\label{eq2.5}
 \limsup_{n\to \infty}|\bk n|^{-1} \ep  \max_{\bk k\le \bk n} T_{\bk k}^2
\le K \ep Y_{\bk 1}^2.
\end{equation}
\end{proposition}

The next theorem gives the  estimate of $\ep M_{\bk n}$ in (\ref{eq2.4}) for  a
non-stationary negatively associated field.
\begin{theorem} \label{theorem2.1}Let $\{ Y_{\bk k};
\bk k \le \bk n \}$ be a negatively associated field of random
variables with $\ep Y_{\bk k}=0$ and $\ep Y_{\bk k}^2<\infty$.
 Denote by $T_{\bk k}= \sum_{\bk i
\le \bk k} Y_{\bk i}$, $M_{\bk n}=\max_{\bk k \le \bk n} |T_{\bk
k}|$. Then there is a constant $C_d$ depending only on $d$ such that
\begin{equation} \label{eqth2.1.1}
\ep M_{\bk n}\le C_d \sqrt{\sum_{\bk k\le \bk n} \ep Y_{\bk k}^2}.
\end{equation}
\end{theorem}

{\bf Proof.} We will prove (\ref{eqth2.1.1}) by induction on $d$. We use the argument of Utev and Peligrad (2003).
For each $\bk n$, define
\begin{equation}\label{eqproofth2.1.1}
a_{\bk n}= \sup_{Y} \left( \ep \max_{\bk m\le \bk n} \left|\sum_{\bk
k\le \bm m}Y_{\bk k}\right|\Big/\left[\sum_{\bk k\le \bk n} \ep
Y_{\bk k}^2\right]^{1/2}\right)
\end{equation}
where the supremum is taken over all fields $Y:=\{Y_{\bk i}\}$ of
square integrable centered negatively associated random variables.

Fix such a random field $\{Y_{\bk i}\}$ and in addition without loss
of generality assume that
$$ \sum_{\bk k\le \bk n} \ep Y_{\bk k}^2=1. $$
Let $M$ be a positive integer that will be specified later. Let
$f(x)=((-M^{-1/2})\vee x)\wedge M^{-1/2}$. For $\bk k\le \bk n$ define:
$ \xi_{\bk k}= f(Y_{\bk k})-\ep f(Y_{\bk k}) $
and
$ \eta_{\bk k}= Y_{\bk k}- \xi_{\bk k}. $ Then both $\{\xi_{\bk n}\}$ and $\{\eta_{\bk n}\}$ are  negatively associated fields.
Since
$$ \sum_{\bk k\le \bk n} \ep |\eta_{\bk k} |
\le 2 M^{1/2} \sum_{\bk k\le \bk n} \ep Y_{\bk k}^2= 2M^{1/2},
 $$
 we
get
$$ \ep \max_{\bk m\le \bk n}\left|\sum_{\bk k\le \bk m}Y_{\bk k} \right|\le
 \ep \max_{\bk m\le \bk n}\left|\sum_{\bk k\le \bk m}\xi_{\bk k} \right|+2 M^{1/2}. $$
To estimate $\ep \max\limits_{\bk m\le \bk n}\left|\sum_{\bk k\le \bk m}\xi_{\bk k} \right|$, we shall use a blocking procedure.

Write $\overline{\bk i}=(i_1,\cdots, i_{d-1})$. Take $m_0=0$ and define the integers $m_k$ recursively by
$$ m_{k}=\min\left\{ m: m>m_{k-1}, \sum_{j=m_{k-1}+1}^m \sum_{\overline{\bk i}\le \overline{\bk n}}\ep
\xi_{\overline{\bk i},j}^2>\frac{1}{M}\right\}. $$
Note that, if we denote by $\ell$ the number of integers produced by this procedure, i.e.: $m_0$, $m_1$, $\ldots$, $m_{\ell-1}$, we have
$$ 1\ge \sum_{k=1}^{\ell -1} \sum_{j=m_{k-1}+1}^{m_k} \sum_{\overline{\bk i}, j}
\ep \xi_{\overline{\bk i}, j}^2>\frac{\ell-1}{M}
\qquad \text{so that } \ell \le M. $$

Write $S_{\overline{\bk m},k}=\sum_{j=m_{k-1}+1}^{m_k}\sum_{\overline{\bk i}\le \overline{\bk m}}\xi_{\overline{\bk i}, j}$
for $1\le k\le \ell-1$ and for convenience, $m_{\ell}=n_d$ that is
$S_{\overline{\bk m},\ell }=\sum_{j=m_{\ell-1}+1}^{n_d}
\sum_{\overline{\bk i}\le \overline{\bk m}}\xi_{\overline{\bk i}, j}$. It is obvious that
\begin{eqnarray}\label{eqproofth2.1.2}
&&\ep \max_{\bk m\le \bk n}\left|\sum_{\bk k\le \bk m}\xi_{\bk k} \right| \nonumber\\
&\le& \ep \max_{1\le k\le \ell}\max_{\overline{\bk m}\le \overline{\bk n}}
\left|\sum_{j=1}^k S_{\overline{\bk m},j} \right|
+\ep \max_{1\le k\le \ell}\max_{\overline{\bk m}\le \overline{\bk n}}
\left(\max_{m_{k-1}<j< m_k}\left|\sum_{t=m_{k-1}+1}^j\sum_{\overline{\bk k}\le \overline{\bk m}}\xi_{\overline{\bk k},t} \right|\right)
\nonumber\\
&=:& I +II.
\end{eqnarray}
We evaluate the two terms in the right hand side of (\ref{eqproofth2.1.2}) separately.

By the induction hypothesis and Cauchy-Schwartz inequality
\begin{eqnarray*}
I&\le& \sum_{k=1}^{\ell}\ep\max_{\overline{\bk m}\le \overline{\bk n}}
\left| S_{\overline{\bk m},k} \right|\le C_{d-1}\sum_{k=1}^{\ell}
\sqrt{ \sum_{\overline{\bk i}\le \overline{\bk n}}\ep\big(\sum_{j=m_{k-1}+1}^{m_k}\xi_{\overline{\bk i},j}\big)^2 }\\
&\le& C_{d-1}\sum_{k=1}^{\ell}
\sqrt{ \sum_{j=m_{k-1}+1}^{m_k}\sum_{\overline{\bk i}\le \overline{\bk n}}\ep \xi_{\overline{\bk i},j}^2 }
\le C_{d-1} M^{1/2}\sqrt{\sum_{\bk i\le \bk n} \ep \xi_{\bk i}^2}\le K_{d-1} M^{1/2}.
\end{eqnarray*}
When $d=1$, the above inequality obviously holds with $C_0=1$, since $\max_{\overline{\bk m}\le \overline{\bk n}}$ does not appear.
To estimate $II$ we notice that
$$ (II)^4\le \sum_{k=1}^{\ell}
\ep \max_{\overline{\bk m}\le \overline{\bk n}}
\max_{m_{k-1}<j< m_k}\left|\sum_{t=m_{k-1}+1}^j\sum_{\overline{\bk k}\le \overline{\bk m}}\xi_{\overline{\bk k},t} \right|^4.
$$
By (\ref{eq2.4}), we obtain
\begin{eqnarray*}
&&\ep \max_{\overline{\bk m}\le \overline{\bk n}}
\max_{m_{k-1}<j< m_k}
\left|\sum_{t=m_{k-1}+1}^j\sum_{\overline{\bk k}\le \overline{\bk m}}\xi_{\overline{\bk k},t} \right|^4\\
&\le& A_{4,d}\left\{ \ep^4 \max_{\overline{\bk m}\le \overline{\bk n}}
\max_{m_{k-1}<j< m_k}\left|\sum_{t=m_{k-1}+1}^j\sum_{\overline{\bk k}\le \overline{\bk m}}\xi_{\overline{\bk k},t} \right|
\right.\\
&& \qquad \left. +\sum_{t=m_{k-1}+1}^{m_k-1}\sum_{\overline{\bk k}\le \overline{\bk n}}\ep \xi_{\overline{\bk k},t}^4
+\left(\sum_{t=m_{k-1}+1}^{m_k-1}\sum_{\overline{\bk k}\le \overline{\bk n}}\ep \xi_{\overline{\bk k},t}^2\right)^2\right\}.
\end{eqnarray*}
By the definition of $m_k$,
$$\sum_{t=m_{k-1}+1}^{m_k-1}\sum_{\overline{\bk k}\le \overline{\bk n}}\ep \xi_{\overline{\bk k},t}^2\le \frac{1}{M},$$
so that by using the notation (\ref{eqproofth2.1.1}) for $a_{\bk n}$ and the definition of $\xi_{\bk k}$ we obtain
$$ \sum_{t=m_{k-1}+1}^{m_k-1}\sum_{\overline{\bk k}\le \overline{\bk n}}\ep \xi_{\overline{\bk k},t}^2
\le \frac{4}{M}\sum_{t=m_{k-1}+1}^{m_k-1}\sum_{\overline{\bk k}\le \overline{\bk n}}\ep \xi_{\overline{\bk k},t}^4\le
\frac{4}{M^2}$$
and
$$\ep^4 \max_{\overline{\bk m}\le \overline{\bk n}}
\max_{m_{k-1}<j< m_k}\left|\sum_{t=m_{k-1}+1}^j\sum_{\overline{\bk k}\le \overline{\bk m}}\xi_{\overline{\bk k},t} \right|
\le a_{\bk n}^4 \left(\sum_{t=m_{k-1}+1}^{m_k-1}\sum_{\overline{\bk k}\le \overline{\bk n}}\ep \xi_{\overline{\bk k},t}^2\right)^2
\le \frac{a_{\bk n}^4}{M^2}. $$
Hence we obtain
$$ (II)^4\le A_{4,d}\left\{ \frac{a_{\bk n}^4}{M}+\frac{4}{M}+\frac{1}{M}\right\}. $$
Now by (\ref{eqproofth2.1.2}) and the estimates for $I$ and $II$ we get
$$ \ep \max_{\bk m\le \bk n}\left|\sum_{\bk k\le \bk m}Y_{\bk k}\right|
\le (A_{4,d}/M)^{1/4} a_{\bk n}+(5A_{4,d}/M)^{1/4} + C_{d-1} M^{1/2}+2M^{1/2}. $$
Therefore, by the definition of $a_{\bk n}$,
$$a_{\bk n}\le (A_{4,d}/M)^{1/4} a_{\bk n}+(5A_{4,d}/M)^{1/4} + C_{d-1} M^{1/2}+2M^{1/2}. $$
Letting $M=[16A_{4,d}]+1$ yields
$$ a_{\bk n}\le 2+2(C_{d-1}+2)M^{1/2}. $$
The proof is now completed.

\bigskip
Combining (\ref{eq2.4}) and (\ref{eqth2.1.1}) we obtain the following result on the moment inequality for the maximum partial sums.
\begin{theorem}\label{theorem2.2}
 Let $p\ge 2$, and let $\{ Y_{\bk k};
\bk k \le \bk n \}$ be a negatively associated field of random
variables with $\ep Y_{\bk k}=0$ and $\ep |Y_{\bk k}|^p<\infty$.
 Then there exists a constant $C_{p,d}$ depending on $p$ and $d$  such that
\begin{equation}\label{eqth2.2.1}
\ep \max_{\bk m\le \bk n}\big| \sum_{\bk k\le \bk m} Y_{\bk k} \big|^p \le C_{p,d} \Big\{ \big( \sum_{\bk k \le \bk n} \ep Y_{\bk k}^2 \big)^{p/2}
+\sum_{\bk k \le \bk n} \ep |Y_{\bk k}|^p \Big\}.
\end{equation}
\end{theorem}

In the case of $d=1$, inequality (\ref{eqth2.2.1}) is first obtained by Shao (2000).
Before that, Su et al (1997) proved that
$$\ep \max_{m\le  n}\big| \sum_{k=1}^m Y_{k} \big|^p \le C_p \Big\{\big(n \max_{k\le n} \ep Y_{k}^2 \big)^{p/2}
+n \sup_{k\le n}\ep |Y_{  k}|^p \Big\}.$$
With the same proof, one can show that (\ref{eqth2.2.1}) of a weakly dependent random field with $\lim_{n\to \infty}\rho^{\ast}_n<1$
(for definition, see Peligrad and Gut (1999)) and a $\rho^{-}$-mixing random field (for definition, see Zhang and Wang (1999) or
Zhang and Wen (2001a)).
%%%%%%%%%%%%%%%%%%%%%%%%%%%%%%%%%%%%%%%%%%%%%%%%%%%%%%%%%%%%%%%%%%%%%%%%

\bigskip
Now, we begin to establish to following Kolmogorov type exponential
inequality for maximum partial sums.

\begin{theorem}\label{prop2.3}
 Let $\{ Y_{\bk k}; \bk k \le \bk n \}$ be a
negatively associated field of random variables with $\ep Y_{\bk
k}=0$ and $ |Y_{\bk k}|\le b$ a.s for some $0<b<\infty$.  Denote
by $T_{\bk k}= \sum_{\bk i \le \bk k} Y_{\bk i}$, $M_{\bk
n}=\max_{\bk k \le \bk n} | T_{\bk k}|$ and $B_{\bk n}^2 =
\sum_{\bk k\le \bk n} \ep Y_{\bk k}^2$. Then for all $x>0$,
\begin{equation}\label{eq2.6}
 \pr ( M_{\bk n} -2\ep M_{\bk n} \ge 20 x)
\le 2^{d+1}\exp \big( - \frac {x^2}{ 2(bx + B_{\bk n}^2) } \big).
\end{equation}
\end{theorem}

\proof  Let $\{ \epsilon_{\bk k}; \bk k \le \bk n \}$ is a field
of i.i.d.r.v.s with $\pr(\epsilon_{\bk k} = \pm 1)=1/2$. Also,
assume that $\{ \epsilon_{\bk k}; \bk k \le \bk n \}$ is
independent of $\{ Y_{\bk k}; \bk k \le \bk n \}$. Denote by
$\widetilde T_{\bk n} =\sum_{\bk k\le \bk n} \epsilon_{\bk
k}Y_{\bk k}$, $ T_{\bk n, 1} =\sum_{\bk k \le \bk n, \epsilon_{\bk
k}=1} Y_{\bk k}$ and $ T_{\bk n, 2} =\sum_{\bk k \le \bk n,
\epsilon_{\bk k}=-1} Y_{\bk k}$. First, we show that
\begin{equation}\label{eq2.7}
 \ep e^{M_{\bk n}}
\le 2^d e^{2 \|M_{\bk n}\|_1}\ep e^{10 |\widetilde T_{\bk n}|} \le
2^{d+1} e^{2 \|M_{\bk n}\|_1}\prod_{\bk k \le \bk n }
   \ep e^{20  \epsilon_{\bk k}Y_{\bk k}}.
\end{equation}

By  the L\'evy inequality, we have
$$\pr\{\widetilde M_{\bk n}\ge x\}
=\ep_Y  \pr_{\epsilon}\{  \widetilde M_{\bk n}\ge x\} \le 2^d
\ep_Y \pr_{\epsilon}\{ |\widetilde T_{\bk n}|\ge x \} \le 2^d
\pr\{ |\widetilde T_{\bk n}|\ge x \}, \quad \forall x\ge 0, $$
where $\ep_Y(\cdot)=\ep[\cdot|\epsilon_{\bk k}, \bk k \in \Bbb
N^d]$, and $\ep _{\epsilon}$, $\pr_{\epsilon}$ etc are defined
similarly.
 Then
$$\ep e^{10 \widetilde M_n}\le 2^d \ep e^{10|\widetilde T_{\bk
n}|}. $$
So,  from (\ref{eq2.3}) it follows that
\begin{eqnarray*}
 \ep e^{M_{\bk n}}  &=& 1+ \ep  M_{\bk n}
 + \sum_{q=2}^{\infty} \frac{ \ep  M_{\bk n}^q }{q! }\\
&\le& 1 + \| M_{\bk n}\|_1 +\sum_{q=2}^{\infty}
\frac{ (  5 \|\widetilde M_{\bk n}\|_q +\|M_{\bk n}\|_1)^q }{q!}\\
&\le& 1 + \| M_{\bk n}\|_1 +\sum_{q=2}^{\infty}
\frac{ \big( 2\big\|5 \widetilde M_{\bk n} +\|M_{\bk n}\|_1\big\|_q\big)^q }{q! }\\
&=& 1 + \| M_{\bk n}\|_1 +\sum_{q=2}^{\infty}
\frac{ \ep( 10 \widetilde M_{\bk n} +2\|M_{\bk n}\|_1)^q }{q! }\\
&\le&1+ \sum_{q=1}^{\infty} \frac{ \ep( 10 \widetilde M_{\bk n} +2
\|M_{\bk n}\|_1)^q }{q! } =\ep\Big\{ 1+ \sum_{q=1}^{\infty}
\frac{ ( 10\widetilde M_{\bk n} +2 \|M_{\bk n}\|_1)^q }{q! } \Big\} \\
&=& \ep e^ { 10 \widetilde M_{\bk n} +2 \|M_{\bk n}\|_1 } \le 2^d
e^{2 \|M_{\bk n}\|_1 }\ep  e^ { 10 |\widetilde T_{\bk n}| }.
\end{eqnarray*}
Note that
$$ \ep e^{10  \widetilde T_{\bk n} } =
\ep e^{10  T_{\bk n,1}- 10  T_{\bk n,2} } \le \frac 12\ep
\big(e^{20  T_{\bk n,1} } + e^{-20 T_{\bk n,2} }\big).
$$
For fixed $\{ \epsilon_{\bk k}; \bk k \le \bk n \}$, we have by
Lemma \ref{lem2.4} or the definition (\ref{eq1.1}),
\begin{eqnarray*}
\ep_Y e^{20 T_{\bk n,1} }
&\le& \prod_{\epsilon_{\bk k}=1} \ep_Y
e^{20 Y_{\bk k}}
=\prod_{\epsilon_{\bk k}=1} \ep_Y e^{20 \epsilon_{\bk k}Y_{\bk k}} \\
&\le& \prod_{\epsilon_{\bk k}=1} \ep_Y e^{20 \epsilon_{\bk
k}Y_{\bk k}} \cdot \prod_{\epsilon_{\bk k}=-1} \ep_Y
e^{20\epsilon_{\bk k}Y_{\bk k}} =\prod_{\bk k\le \bk n } \ep_Y
e^{20 \epsilon_{\bk k}Y_{\bk k}},
\end{eqnarray*}
since $\ep_Y e^{20 \epsilon_{\bk k}Y_{\bk k}} \ge e^{\ep_Y (20
\epsilon_{\bk k}Y_{\bk k}) }=1$. It follows that
\begin{eqnarray*}
\ep e^{20 T_{\bk n,1}}
&=& \ep_{\epsilon}\ep_Y e^{20 T_{\bk n,1}}
\le \ep_{\epsilon}\big(\prod_{\bk k\le \bk n }
   \ep_Y e^{20 \epsilon_{\bk k}Y_{\bk k}}\big) \\
&=&\prod_{\bk k\le \bk n } \ep_{\epsilon}
   \ep_Y e^{20 \epsilon_{\bk k}Y_{\bk k}}
=\prod_{\bk k\le \bk n } \ep e^{20 \epsilon_{\bk k}Y_{\bk k}}.
\end{eqnarray*}
Similarly,
$$ \ep e^{-20 T_{\bk n,2}}
\le \prod_{\bk k\le \bk n } \ep e^{-20 \epsilon_{\bk k}Y_{\bk k}}
=\prod_{\bk k\le \bk n } \ep e^{20 \epsilon_{\bk k}Y_{\bk k}}. $$
It follows that
$$ \ep e^{10 \widetilde T_{\bk n}}
\le \prod_{\bk k\le \bk n } \ep e^{20 \epsilon_{\bk k}Y_{\bk k}}.
$$
Similarly,
$$\ep e^{-10 \widetilde T_{\bk n}}
\le \prod_{\bk k\le \bk n } \ep e^{20 \epsilon_{\bk k}Y_{\bk k}}.
$$
(\ref{eq2.7}) is proved.
 Now, from (\ref{eq2.7}) it follows that for any $x>0$ and $t>0$,
$$ \pr(M_{\bk n}-2\ep M_{\bk n} \ge 20 x)
\le e^{-20 tx - 2t \| M_{\bk n}\|_1 }\ep e^{t M_{\bk n} } \le
2^{d+1} e^{-20 tx} \prod_{\bk k\le \bk n }
   \ep e^{20 t \epsilon_{\bk k}Y_{\bk k}}.  $$
Since
\begin{eqnarray*}
\ep e^{20 t \epsilon_{\bk k} Y_{\bk k} }
&=& \ep\Big\{1+ 20 t
\epsilon_{\bk k} Y_{\bk k} +\frac{ e^{20 t \epsilon_{\bk k}Y_{\bk
k} }-1
   -20 t \epsilon_{\bk k}Y_{\bk k} }{ (\epsilon_{\bk k} Y_{\bk k} )^2}
(\epsilon_{\bk k}Y_{\bk k} )^2\Big\} \\
&\le& 1+ (e^{20tb}-1- 20 tb) b^{-2} \ep Y_{\bk k}^2 \\
&\le&\exp \big\{(e^{20tb}-1- 20 tb) b^{-2} \ep Y_{\bk k}^2\big\},
\end{eqnarray*}
it follows that
$$\pr(M_{\bk n}- 2\ep M_{\bk n} \ge 20 x)
\le 2^{d+1} e^{-20 tx} \exp \big\{(e^{20tb}-1- 20 tb) b^{-2}
B_{\bk n}^2\big\}. $$ Letting $20 t =\frac 1b\log
(1+\frac{xb}{B_{\bk n}^2})$ yields
\begin{eqnarray*}
\pr(M_{\bk n}-2\ep M_{\bk n} \ge 20 x)
&\le& 2^{d+1} \exp\Big\{
\frac xb -(\frac xb+ \frac{ B_{\bk n}^2}{b^2})
  \log(1+\frac{xb}{ B_{\bk n}^2 }) \Big\} \\
&\le& 2^{d+1}\exp\Big\{ - \frac{x^2}{2(bx+ B_{\bk n}^2)}\Big\}.
\end{eqnarray*}

%%%%%%%%%%%%%%%%%%%%%%%%%%%%%%%%%%%%%%%%%%%%%%%%%%%%%%%%%%%%%%
\section{Proofs of the laws of iterated logarithm}
\setcounter{equation}{0}

We need two more lemmas.
\begin{lemma}\label{lem2.8}
Let $d\ge 2$ and let $\{ Y_{\bk k}; \bk
k \in \Bbb N^d \}$ be a negatively associated field of identically
distributed random variables with
$$\ep Y_{\bk 1}=0 \; \text {and }\;
\ep Y_{\bk 1}^2\log^{d-1}(|Y_{\bk 1}|)/\log\log(|Y_{\bk
1}|)<\infty. $$
 Denote by $T_{\bk n}= \sum_{\bk k\le \bk n} Y_{\bk
k}$ and $M_{\bk n}=\max_{\bk k\le \bk n}|T_{\bk k}|$. Then
\begin{eqnarray}\label{eq2.8}
\limsup_{\bk n\to \infty} \frac {M_{\bk n}}{(2d |\bk n| \log
\log |\bk n|)^{1/2}}
\le 20 (\ep Y_{\bk 1}^2)^{1/2} \quad a.s.
\end{eqnarray}
\end{lemma}

\proof Let $0<\epsilon<1$ be an arbitrary but fixed number. Let
$b_m= \frac{\epsilon}{40}(\ep Y_{\bk 1}^2)^{1/2}
 (m/\log\log m)^{1/2}$, $f_{\bk k}(x)=(- b_{|\bk k|})\vee (x \wedge b_{|\bk k|})$,
$g_{\bk k}(x)= x- f_{\bk k}(x)$. Define
\begin{eqnarray*}
&& \overline Y_{\bk k}= f_{\bk k}(Y_{\bk k})- \ep f_{\bk k}(Y_{\bk
k}), \quad
\widehat Y_{\bk k}= g_{\bk k}(Y_{\bk k}), \\
&& \overline T_{\bk n}=\sum_{\bk k\le \bk n } \overline Y_{\bk k},
\quad
 \overline M_{\bk n}= \max_{\bk k \le \bk n} | \overline T_{\bk k}|,
\quad   \widehat T_{\bk n}=\sum_{\bk k\le \bk n} \widehat Y_{\bk
k}.
\end{eqnarray*}
First, we show that
\begin{equation}\label{eq2.9}
 \frac{
\widehat T_{\bk n} - \ep \widehat T_{\bk n}}{ (2|\bk n| \log \log
|\bk n|)^{1/2}} \to 0 \quad a.s. \quad \text{ as } \bk n \to
\infty. \end{equation}
 Since
\begin{eqnarray}\label{eq2.10}
&& \frac {\sum_{\bk k\le \bk n}\ep |\widehat
Y_{\bk k}|}{ (2|\bk n| \log \log |\bk n|)^{1/2}} \le  \frac
{\sum_{\bk k\le \bk n}\ep | Y_{\bk k}| I\{ |Y_{\bk k}| \ge b_{|\bk
k|}\} }{
      (2|\bk n| \log \log |\bk n|)^{1/2}} \non\\
&& \quad  \le C \frac {\sum_{\bk k\le \bk n} (\frac{ \log\log |\bk
k|}{|\bk k|})^{1/2}
        \ep  Y_{\bk 1}^2 I\{ |Y_{\bk 1}| \ge b_{|\bk k|} \} }{
        (2|\bk n| \log \log |\bk n|)^{1/2}} \non\\
&& \quad  \le C \frac {\sum_{\bk k\le \bk n} |\bk k|^{-1/2}
        \ep  Y_{\bk 1}^2 I\{ |Y_{\bk 1}| \ge b_{|\bk k|} \} }{
        |\bk n|^{1/2}} \non\\
&& \quad  = o(1) \frac {\sum_{\bk k\le \bk n} |\bk k|^{-1/2}
         }{
        |\bk n|^{1/2}}
\to 0 \quad \text{ as } |\bk n|\to \infty,
\end{eqnarray}
it is enough to show that
\begin{equation}\label{eq2.11}
 \frac
{\sum_{\bk k\le \bk n}|\widehat Y_{\bk k}|}{ (2|\bk n| \log \log
|\bk n|)^{1/2}} \to 0 \quad a.s.
\end{equation}
 For $\bk
n=(n_1,\cdots, n_d)\in \Bbb N^d$, let $I(\bk n)=\{ \bk
k=(k_1,\cdots, k_d): 2^{n_i-1}\le k_i \le 2^{n_i}-1, i=1,\cdots
d\} $. (\ref{eq2.11}) will be true if we have
\begin{equation}\label{eq2.12}
\frac{ \sum_{\bk k\in I(\bk n)} \widehat Y_{\bk k}^+ }{
    (2^{\|\bk n\|} \log\log 2^{\|\bk n\|})^{1/2} }
\to 0 \quad a.s.
\end{equation}
and
\begin{equation}\label{eq2.13}
 \frac{ \sum_{\bk k\in I(\bk n)} \widehat Y_{\bk k}^- }{
    (2^{\|\bk n\|} \log\log 2^{\|\bk n\|})^{1/2} }
\to 0 \quad a.s.
\end{equation}
We show (\ref{eq2.12}) only since (\ref{eq2.13}) can be showed
similarly. Let $\alpha_m :=\alpha(m)=(2m \log \log m)^{1/2}$ and
$Z_{\bk k}=(\widehat Y_{\bk k}^+)\wedge \alpha_{|\bk k|}$. It is
easily seen that $Z_{\bk k}=0$ if $Y_{\bk k}\le b_{|\bk k|}$,
$Z_{\bk k}=Y_{\bk k}-b_{|\bk k|}$ if $b_{|\bk k|}\le Y_{\bk k}\le
b_{|\bk k|}+\alpha_{|\bk k|}$ and $Z_{\bk k}=\alpha_{|\bk k|}$ if
$Y_{\bk k}\ge b_{|\bk k|}+\alpha_{|\bk k|}$. Also, $\{Z_{\bk k};
k\in \Bbb N^d\}$ is a negatively associated field of random
variables. Obviously,
\begin{eqnarray*}
\sum_{\bk n} \pr(\widehat Y_{\bk n}^+ \ne Z_{\bk n})
&\le&\sum_{\bk n}\pr(Y_{\bk n}\ge \alpha_{|\bk n|})
\le  C \sum_{m=1}^{\infty} (\log m)^{d-1}\pr(Y_{\bk 1}\ge \alpha_m)\\
&\le& C \ep Y_{\bk 1}^2\log^{d-1}(|Y_{\bk 1}|)/ \log\log(|Y_{\bk
1}|)<\infty.
\end{eqnarray*}
Also by (\ref{eq2.10}), using the notation $2^{\bk
n}=(2^{n_1},\cdots, 2^{n_d})$ for $\bk n= (n_1,\cdots, n_d)\in
\Bbb N^d$, one has that,
$$ \frac{\big|\sum_{\bk k\in I(\bk n)} \ep Z_{\bk k} \big| }{
    (2^{\|\bk n\|} \log\log 2^{\|\bk n\|})^{1/2} }
\le \frac{ \sum_{\bk k\le 2^{\bk n} } \ep|\widehat  Y_{\bk k}\big| }{
    (2^{\|\bk n\|} \log\log 2^{\|\bk n\|})^{1/2} }
\to 0. $$ So, (\ref{eq2.12}) is equivalent to
\begin{equation}\label{eq2.14} \frac{ \sum_{\bk
k\in I(\bk n)} (Z_{\bk k} -\ep Z_{\bk k}) }{
    (2^{\|\bk n\|} \log\log 2^{\|\bk n\|})^{1/2} }
\to 0 \quad a.s.
\end{equation}
Let
$$\Lambda(\bk n) =\sum_{\bk k \in I(\bk n) }
\frac{\ep Y_{\bk 1}^2 I\{ b_{|\bk k|}< |Y_{\bk 1}|\le 2 \alpha_{|\bk k|}\} }{\alpha_{|\bk k|}^2}.
$$
Note that $b_{|\bk k|} < \alpha_{|\bk k|}$ if $|\bk k|$ is large
enough. From Lemma \ref{lem2.1}, it follows that for $\|\bk n\|$
large enough and any $\delta>0$,
\begin{eqnarray*}
&& \pr \big(\big| \sum_{\bk k\in I(\bk n)} (Z_{\bk k} -\ep Z_{\bk
k}) \big| \ge \delta (2^{\|\bk n\|} \log \log 2^{\|\bk n\|})^{1/2}
\big)\\
&\le& C \frac{\ep \big| \sum_{\bk k\in I(\bk n)}
    (Z_{\bk k} -\ep Z_{\bk k}) \big|^4}{
     (2^{\|\bk n\|} \log \log 2^{\|\bk n\|})^2 } \\
&\le& C (\alpha(2^{\|\bk n\|}) )^{-4} \Big\{ \big(\sum_{\bk k \in
I(\bk n) }\ep |Z_{\bk k}|^2\big)^2
   +\sum_{\bk k \in I(\bk n) }\ep |Z_{\bk k}|^4 \Big\} \\
&\le&C (\alpha(2^{\|\bk n\|}) )^{-4} \Big\{ \big(\sum_{\bk k \in
I(\bk n)    } \ep  Y_{\bk 1}^2 I\{   b_{|\bk k|} <|Y_{\bk 1}|
   \le \alpha_{|\bk k|}+b_{|\bk k|} \} \big)^2 \\
&&\qquad   +\sum_{\bk k \in I(\bk n) }\ep Y_{\bk 1}^4 I\{
b_{|\bk k|} <|Y_{\bk 1}|\le \alpha_{|\bk k|}+b_{|\bk k|}\}
 \Big\} \\
&&+C (\alpha(2^{\|\bk n\|}) )^{-4} \Big\{ \big(\sum_{\bk k \in
I(\bk n)    }
\alpha_{|\bk k|}^2 \pr(  |Y_{\bk 1}|\ge \alpha_{|\bk k|}+b_{|\bk k|} ) \big)^2  \\
&& \qquad  +\sum_{\bk k \in I(\bk n) }\alpha_{|\bk k|}^4
\pr(|Y_{\bk 1}|\ge \alpha_{|\bk k|}+b_{|\bk k|})
 \Big\} \\
&\le& C\Big\{ \Lambda^2(\bk n) + \sum_{\bk k \in I(\bk n)}
\frac{\ep Y_{\bk 1}^4 I\{ b_{|\bk k|}< |Y_{\bk 1}|\le
2\alpha_{|\bk k|}\} }{ \alpha_{|\bk k|}^4}
\Big\} \\
&&+C\Big\{ \big( \sum_{\bk k \in I(\bk n)} \pr (|Y_{\bk 1}|\ge
\alpha_{|\bk k|})\big)^2 +\sum_{\bk k \in I(\bk n)} \pr (|Y_{\bk
1}|\ge \alpha_{|\bk k|}) \Big\}.
\end{eqnarray*}
Obviously,
\begin{eqnarray*}
&&\sum_{\bk n}  \Big\{ \big( \sum_{\bk k \in I(\bk n)} \pr
(|Y_{\bk 1}|\ge \alpha_{|\bk k|})\big)^2
+\sum_{\bk k \in I(\bk n)} \pr (|Y_{\bk 1}|\ge \alpha_{|\bk k|}) \Big\} \\
&\le& C \sum_{\bk n}\pr (|Y_{\bk 1}|\ge \alpha_{|\bk n|}) \le C
\ep Y_{\bk 1}^2\log^{d-1}(|Y_{\bk 1}|)/ \log \log (|Y_{\bk
1}|)<\infty.
\end{eqnarray*}
Also,
\begin{eqnarray*}
&&\sum_{\bk n}  \sum_{\bk k \in I(\bk n)} \frac{\ep Y_{\bk 1}^4
I\{ b_{|\bk k|}< |Y_{\bk 1}|\le 2\alpha_{|\bk k|}\} }{
\alpha_{|\bk k|}^4}
\le \sum_{\bk n}\frac{\ep Y_{\bk 1}^4 I\{ b_{|\bk n|}< |Y_{\bk 1}|\le 2\alpha_{|\bk n|}\} }{ \alpha_{|\bk n|}^4} \\
&\le&\sum_{m=1}^{\infty} (\log m)^{d-1}
  \frac{\ep Y_{\bk 1}^4 I\{  |Y_{\bk 1}|\le 2\alpha_m \} }{ \alpha_m^4} \\
&\le&C \sum_{m=1}^{\infty} (\log m)^{d-1} \sum_{k=1}^{m}
  \frac{\ep Y_{\bk 1}^4 I\{ 2\alpha_{k-1}<  |Y_{\bk 1}|\le 2\alpha_k \} }{ \alpha_m^4} \\
&\le& C\sum_{k=1}^{\infty}  \sum_{m=k}^{\infty} \frac{(\log
m)^{d-1}}{(m\log\log m)^2}
  \ep Y_{\bk 1}^4 I\{ 2\alpha_{k-1}<  |Y_{\bk 1}|\le 2\alpha_k \} \\
&\le& C \sum_{k=1}^{\infty}\frac{(\log k)^{d-1}}{ k(\log\log k)^2}
(k\log\log k)
 \ep Y_{\bk 1}^2 I\{ 2\alpha_{k-1}<  |Y_{\bk 1}|\le 2\alpha_k \} \\
&\le&C  \ep Y_{\bk 1}^2 \log^{d-1}(|Y_{\bk 1}|)/\log \log (|Y_{\bk
1}|)<\infty,
\end{eqnarray*}
and, similarly to (3.8) of Li and Wu (1989) we have
$$ \sum_{\bk n} \Lambda^2(\bk n)<\infty. $$
It follows that for arbitrary $\delta>0$,
$$ \sum_{\bk n}
\pr \big(  \big| \sum_{\bk k\in I(\bk n)} (Z_{\bk k} -\ep Z_{\bk
k}) \big| \ge \delta (2^{\|\bk n\|} \log \log 2^{\|\bk
n\|})^{1/2}\big) <\infty, $$ which implies (\ref{eq2.14}) by the
Borel-Cantelli lemma. Thus (\ref{eq2.9}) holds.

Now,  by applying (\ref{eq2.6}) to $x= (1+2 \epsilon) (\ep Y_{\bk
1}^2)^{1/2}( 2 d |\bk n| \log\log |\bk n|)^{1/2}$ and $b=2b_{|\bk
n|}$, it follows that
\begin{eqnarray*}
&& \pr \big(  \overline M_{\bk n} -2 \ep  \overline M_{\bk n}\ge
20(1+2 \epsilon) (\ep Y_{\bk 1}^2)^{1/2}( 2 d |\bk n| \log\log |\bk n|)^{1/2} \big)\\
&\le& 2^{d+1} \exp\{ -(1+\epsilon) d \log\log |\bk n|\} \le
2^{d+1}(\log |\bk n|)^{-(1+\epsilon) d}.
\end{eqnarray*}
For $\theta>1$ and $\bk m\in \Bbb N^d$, let $\bk N_{\bk
m}=([\theta^{m_1}], \cdots, [\theta^{m_d}])$. It follows that
\begin{eqnarray*}
&&\sum_{\bk m} \pr \big(  \overline M_{\bk N_{\bk m} } -2 \ep
\overline M_{\bk N_{\bk m} } \ge
20(1+2 \epsilon) (\ep Y_{\bk 1}^2)^{1/2}( 2 d |\bk N_{\bk m}| \log\log |\bk N_{\bk m}|)^{1/2} \big)\\
&\le& C \sum_{\bk m} \| \bk m\|^{-(1+\epsilon)d} \le C
\sum_{i=1}^{\infty} i^{d-1} i^{-(1+\epsilon)d}<\infty.
\end{eqnarray*}
Notice $\ep \overline M_{\bk N_{\bk m} }\le C |\bk N_{\bk m}|^{1/2}$ by Theorem \ref{theorem2.1}.
From the Borel-Cantelli lemma, it follows that
\begin{eqnarray*}
&&\limsup_{ \bk m\to \infty} \frac{ \overline M_{\bk N_{\bk m} }
}{
   (2d  |\bk N_{\bk m}| \log\log |\bk N_{\bk m}|)^{1/2} } \non\\
&\le& 20 (1+2\epsilon)( \ep Y_{\bk 1}^2)^{1/2} + 2 \limsup_{ \bk
m\to \infty} \frac{ \ep \overline M_{\bk N_{\bk m}} }{
   (2d  |\bk N_{\bk m}| \log\log |\bk N_{\bk m}|)^{1/2} }\non\\
 &=& 20 (1+2\epsilon)( \ep Y_{\bk 1}^2)^{1/2}.
\end{eqnarray*}
So,  we have
\begin{eqnarray}\label{eq2.17}
&& \limsup_{ \bk n \to \infty} \frac{  \overline M_{\bk n}}{
   (2d  |\bk  n| \log\log |\bk n| )^{1/2} } \non\\
&\le& \limsup_{ \bk n \to \infty}\max_{\bk N_{\bk m-\bk 1}\le \bk
n \le \bk N_{\bk m}} \frac{ \overline M_{N_{\bk m} } }{
   (2d  |\bk N_{\bk m-\bk 1}| \log\log |\bk N_{\bk m-\bk 1}|)^{1/2} } \non\\
&\le& 20 \theta^{d/2}(1+2\epsilon)( \ep Y_{\bk 1}^2)^{1/2} \quad
a.s.
\end{eqnarray}
 Finally, from (\ref{eq2.9}) and
(\ref{eq2.17}) it follows that (\ref{eq2.8}) holds.

\begin{lemma}\label{lem2.9}
 Let $\{ X_{\bk k}; \bk k \in \Bbb N^d \}$ be a
negatively associated field of bounded random variables with $\ep
X_{\bk k}=0$ for all $\bk k\in \Bbb N^d$. Denote by $S_{\bk n}=
\sum_{\bk k\le \bk n} X_{\bk k}$. Then  for any $\delta>0$,
\begin{equation}\label{eq2.18}
 \limsup_{\bk n\to \infty} \frac{\max_{\bk k\le \delta \bk
n}| S_{\bk n+\bk k}- S_{\bk n}|}{ (2d |\bk n|\log \log |\bk
n|)^{1/2}} \le 80d (\delta (1+2\delta)^d)^{1/2}\sup_{\bk k}(\ep
X_{\bk k}^2)^{1/2}
\end{equation}
\end{lemma}

\proof Assume that $|X_{\bk k}|\le b $ a.s. with $0<b<\infty$.
Denote $\beta=\sup_{\bk k}(\ep X_{\bk k}^2)^{1/2}$. From  Theorem \ref{theorem2.2}  it
follows that
\begin{equation}\label{eq2.19}
 \ep
\max_{\bk m\le \bk n} (\sum_{\bk k\le \bk m} X_{\bk k})^2 \le C |\bk n| \beta^2.
\end{equation}
 For fixed $\bk n$ and $\bk m\le \bk n$, let
$Y_{\bk k}=0$ for $\bk k\le \bk m$ and $Y_{\bk k}= X_{\bk k}$
otherwise, and let $T_{\bk k}=\sum_{\bk i\le \bk k} Y_{\bk i}$ for
$\bk k\le \bk m+\delta \bk n$, $B_{\bk m+[\delta \bk n]}^2
=\sum_{\bk k\le \bk m+[\delta \bk n]} \ep Y_{\bk k}^2$. Then for
$\bk m\le \bk n$,
$$
B_{\bk m+[\delta \bk n]}^2 \le  \big( |\bk m+[\delta \bk n] |-
|\bk m|\big) \beta^2 \le \delta d (1+\delta)^d |\bk n| \beta^2.
$$
So by Theorem \ref{prop2.3},
\begin{eqnarray}\label{eq2.20}
&&\pr\big(     \max_{\bk k\le \delta \bk n}| S_{\bk m+\bk k}-
S_{\bk m}|
 - 2 \ep \max_{\bk k\le \delta \bk n}| S_{\bk m+\bk k}- S_{\bk m}|
  \ge 20 x \big) \non\\
&& =\pr\big(     \max_{\bk k\le \bk m +[\delta \bk n] }|T_{\bk k}|
- \ep \max_{\bk k\le \bk m +[\delta \bk n] }|T_{\bk k}| \ge 20 x \big) \non\\
&& \le 2^{d+1} \exp\Big\{ -\frac {x^2}{ 2( bx + B_{\bk m+[\delta \bk n]}^2) } \Big\} \non\\
&& \le 2^{d+1}\exp\Big\{ -\frac {x^2}{ 2( bx + \delta
d(1+\delta)^d |\bk n| \beta^2 ) } \Big\}. \end{eqnarray} Note that
by (\ref{eq2.19}) we have
\begin{eqnarray*}
&& \max_{\bk m\le n}  \ep \max_{\bk k\le \delta \bk n} | S_{\bk
m+\bk k}- S_{\bk m}| \le  2 \ep \max_{\bk k\le \bk n+[\delta \bk
n]}
   | S_{\bk k}|
\le  2 \big(\ep \max_{\bk k\le \bk n+[\delta \bk n]} S_{\bk k}^2\big)^{1/2} \\
&\le&  K \big(|\bk n+[\delta \bk n]|\beta^2\big)^{1/2}
 \le K (1+\delta)^{d/2} |\bk n|^{1/2} \beta.
\end{eqnarray*}
From (\ref{eq2.20}) it follows that for $|\bk n|$ large enough
\begin{eqnarray*}
&&\max_{\bk m\le \bk n} \pr\Big(   \max_{\bk k\le \delta \bk n}|
S_{\bk m+\bk k}- S_{\bk m}|
 \ge 20d  (1+2\epsilon)
  \big( 2\delta (1+\delta)^d  d^2 \beta^2
   |\bk n|\log \log |\bk n|\big)^{1/2}\Big) \non\\
&&\qquad\le 2^{d+1} \exp\big\{ -(1+\epsilon) d \log \log |\bk n|
\big\}.
\end{eqnarray*}
Now, let $I_{\bk p}=\{ \bk p+\bk k: \bk k\le [\delta \bk n]\}$.
Then there are at most
$[(\delta^{-1}+1)^d]+1$ such $I_{\bk p}$s whose union covers
$\{\bk k:\bk k\le \bk n\}$.
It follows that
\begin{eqnarray*}
&&\pr \Big(   \max_{\bk m\le \bk n} \max_{\bk k\le \delta \bk n}
  |S_{\bk m+\bk k}-S_{\bk m}|
\ge (1+2\epsilon) 80d  \big( \delta (1+2\delta)^d
\beta^2 |\bk n| \log\log |\bk n|\big)^{1/2} \Big) \\
&&\le  \Big( (\frac 1{\delta}+1)^d +1\Big) \max_{\bk p} \pr\Big(
\max_{\bk m \in I_{\bk p} }
\max_{\bk k \le \delta \bk n} |S_{\bk m+\bk k}-S_{\bk m}| \\
&& \qquad \ge 40d(1+2\epsilon)
  \big( 4\delta (1+2\delta)^d   \beta^2
  |\bk n| \log\log |\bk n|\big)^{1/2} \Big) \\
&&\le 4(\frac 1{\delta}+1)^d \max_{\bk m\le \bk n}
\pr\Big(  \max_{\bk k \le 2\delta \bk n } |S_{\bk m+\bk k}-S_{\bk m}| \\
&& \qquad  \ge 20d (1+2\epsilon) \big( 4\delta (1+2\delta)^d
\beta^2
  |\bk n| \log\log |\bk n|\big)^{1/2} \Big) \\
&&\le 4(\frac 1{\delta}+1)^d 2^{d+1}\exp\big(-(1+\epsilon) d \log
\log |\bk n|\big).
\end{eqnarray*}
If we choose $\bk n_{\bk p}=([\theta^{p_1}], \cdots,
[\theta^{p_d}])$, then the sum of the above probability is finite.
And then by the Borel-Cantelli lemma,
$$ \limsup_{\bk p\to \infty}
\frac{ \max_{\bk m\le \bk n_{\bk p}} \max_{\bk k\le \delta \bk n_{\bk p}}
  |S_{\bk m+\bk k}-S_{\bk m}|}{ (2 d|\bk n_{\bk p}| \log\log |\bk n_{\bk p}|)^{1/2} }
\le  80d \big( \delta (1+2\delta)^d\big)^{1/2} \beta\quad a.s., $$
which implies (\ref{eq2.18}) easily.

%%%%%%%%%%%%%%%%%%%%%%%

\bigskip
Now, we turn to the

{\noindent \bf Proof of Theorem \ref{th2}:}
We can assume $\sigma>0$, for otherwise, we can consider the field
$\{X_{\bk n}+\epsilon Z_{\bk n}; \bk n\in \Bbb N^d\}$ instead,
where $\{Z_{\bk n}; \bk n\in \Bbb N^d\}$ is a field of i.i.d.
standard normal random variables and $\epsilon>0$ is an arbitrary
number. First we show that
\begin{equation}\label{upbound}
\limsup_{ \bk n\to \infty} \frac{ |S_{\bk n}| }{ (2d |\bk n| \log
\log |\bk n|)^{1/2} } \le \sigma \quad a.s.
\end{equation}

For $b>0$, let
$g_b(x)=(-b)\vee(x\wedge b)$ and $ h_b(x)=x-g_b(x)$. Then $g_b(x)$
and $h_b(x)$ are both non-decreasing functions of $x$. Let $
\overline X_{\bk k}= g_b(X_{\bk k} )-\ep g_b(X_{\bk k})$,
 $\widehat X_{\bk k} = h_b(X_{\bk k})- \ep h_b(X_{\bk k})$,
  $\overline S_{\bk n}=\sum_{\bk k\le \bk n}\overline X_{\bk k}$
 and $\overline M_{\bk n}=\max_{\bk k\le \bk n} |\overline S_{\bk
 k}|$. And define $\widehat S_{\bk k}$, $\widehat M_{\bk n}$
 similarly.  Then $ X_{\bk k}= \overline X_{\bk k}
+\widehat X_{\bk k}$ and,
 $\{ \overline X_{\bk k}; \bk k\in \Bbb N^d \}$ and $\{\widehat X_{\bk k}; \bk k \in \Bbb N^d\}$
are both  negatively associated fields of
identically distributed random variables with $\ep \overline X_{\bk k} =\ep \widehat
X_{\bk k}= 0$ and $|\overline X_{\bk k}|\le 2b $. Then  by Lemma
\ref{lem2.8},
\begin{eqnarray}\label{eqxhat}
&&\limsup_{\bk n\to \infty} \frac{|\widehat S_{\bk n}|}{ (2d |\bk n| \log\log |\bk n|)^{1/2} }
\nonumber\\
&\le& 20(\ep \widehat X_{\bk 1}^2)^{1/2}
\le  20 \big( \ep X_{\bk 1}^2I\{ |X_{\bk 1}|\ge b\}\big)^{1/2}
 \to 0 \quad a.s. \quad  \text{ as } b\to \infty.
\end{eqnarray}
So it suffices to show that for any $\epsilon>0$,
\begin{equation}\label{eqadP1.2}
\limsup_{ \bk n\to \infty} \frac{ |\overline{S}_{\bk n}| }{ (2d |\bk n| \log
\log |\bk n|)^{1/2} } \le (1+\epsilon)^2 \sigma \quad a.s.
\end{equation}
if $b$ is large enogh.

Let $I_m=\{\bk k: 1\le k_i \le m, i=1,\cdots, d\}$ and $Y_{\bk
k}=\sum_{\bk i\in I_m}\overline{X}_{(\bk k -\bk 1)m+\bk i} $.
Since
$$ \ep(\sum_{\bk i\in I_m}X_{(\bk k -\bk 1)m+\bk i})^2=\ep S_{ \bk 1m }^2 /m^d \to \sigma^2 $$
and
$$\ep\Big(Y_{\bk k}-\sum_{\bk i\in I_m}X_{(\bk k -\bk 1)m+\bk i}\Big)^2/m^d\le 2
 \ep X_{\bk 1}^2I\{ |X_{\bk 1}|\ge b\}
\to 0  \quad \text{ as } b\to \infty,$$
we can choose $b$ and $m$ large
enough such that
$$ \sup_{\bk k} \ep
Y_{\bk  k}^2  \le (1+\epsilon)^2 m^d \sigma^2.   $$
 For $\theta>1$,
let $\bk N_{\bk k}=[\theta^{\bk k}]=:
  ( [\theta^{k_1}], \cdots, [\theta^{k_d}] )$.
From (\ref{eq2.1}), it follows that
$$\pr(  |\overline{S}_{\bk N_{\bk k}m} | \ge x)=\pr(|\sum_{\bk i\le \bk
N_{\bk k}}Y_{\bk i}|\ge x)\le  2 \exp
  \Big\{ -\frac {x^2}{ 2( 2m^dbx +\sum_{\bk i\le \bk  N_{\bk k}} \ep Y_{\bk i}^2) }\Big\}.
$$
Let $ x= (1+\epsilon)^2 \sigma
  (2 d m^d |\bk N_{\bk k}| \log\log |\bk N_{\bk k}|)^{1/2}$.
It follows that
\begin{eqnarray*}
&&\sum_{\bk k}   \pr\big( | \overline{S}_{\bk N_{\bk k} m  }|
 \ge  (1+\epsilon)^2 \sigma
    (2 d m^d |\bk N_{\bk k}|  \log\log |\bk N_{\bk k}|)^{1/2} \big)\\
&\le &\sum_{\bk k}   \pr\big( |\overline{S}_{\bk N_{\bk k} m  }|
 \ge  (1+\epsilon)
    (2 d \sum_{\bk i\le \bk N_{\bk k}} \ep Y^2_{\bk i}  \log\log |\bk N_{\bk k}|)^{1/2} \big) \\
&\le& C\sum_{\bk k} \exp\big\{ -(1+\epsilon) d \log \log |\bk
N_{\bk k}| \big \} \\
&\le& C \sum_{\bk k} \|\bk k\|^{-(1+\epsilon) d}
\le C \sum_{i=1}^{\infty} i^{d-1} i^{-(1+\epsilon)d} <\infty.
\end{eqnarray*}
From the Borel-Cantelli lemma, it follows that
\begin{eqnarray*}
&&\limsup_{\bk k\to \infty} \frac{ |\overline{S}_{\bk N_{\bk k} m} |}{
  (2d |\bk N_{\bk k}m | \log\log |\bk N_{\bk k} m|)^{1/2} } \\
&=& \limsup_{\bk k\to \infty} \frac{ |\overline{S}_{\bk
N_{\bk k} m} |}{
  (2d m^d |\bk N_{\bk k} | \log\log |\bk N_{\bk k} |)^{1/2} }
\le   (1+\epsilon)^2\sigma \quad
a.s.
\end{eqnarray*}
 So, from Lemma \ref{lem2.9} it follows that
\begin{eqnarray*}
&&\limsup_{ \bk n\to \infty}
 \frac{ | \overline{S}_{\bk n}|}{(2d |\bk n|\log \log |\bk n|)^{1/2}} \\
&\le& \limsup_{ \bk k\to \infty} \max_{\bk N_{\bk k}m  \le \bk n
\le \bk N_{\bk k+\bk 1}m }
 \frac{ |\overline{S}_{\bk n}|}{(2d |\bk N_{\bk k}m | \log \log |\bk N_{\bk k}m| )^{1/2}} \\
&\le& \limsup_{ \bk k\to \infty}
 \frac{ | \overline{S}_{\bk N_{\bk k}m }|}{(2d |\bk N_{\bk k}m | \log \log |\bk N_{\bk k} m|)^{1/2}} \\
&& \quad +\limsup_{ \bk k\to \infty} \max_{\bk N_{\bk k}m  \le \bk
n \le \bk N_{\bk k+\bk 1}m }
 \frac{ | \overline{S}_{\bk n}-  \overline{S}_{\bk N_{\bk k}m }| }{
    (2d |\bk N_{\bk k}m | \log \log |\bk N_{\bk k}m |)^{1/2}} \\
&\le&(1+\epsilon)\sigma+ 80d \big( (\theta-1) (1+2(\theta-1))^d
\big)^{1/2} (\ep \overline{X}^2)^{1/2} \quad a.s.
\end{eqnarray*}
Letting $\theta\to 1$ completes the proof of
(\ref{eqadP1.2}).

Next we  show that
$$
\limsup_{ \bk n\to \infty} \frac{ S_{\bk n} }{ (2d |\bk n| \log \log
|\bk n|)^{1/2} } \ge  \sigma \quad a.s.
$$
Due to (\ref{eqxhat}), it suffices to show that for any
$0<\epsilon<1/9$,
\begin{equation}\label{eq2.21}
\limsup_{ \bk n\to \infty} \frac{ \overline{S}_{\bk n} }{ (2d |\bk
n| \log \log |\bk n|)^{1/2} } \ge (1-9 \epsilon) \sigma \quad a.s.
\end{equation}
for $b$ large enough.

 For $k \in \Bbb N$, let $m_k=[ 2^{k^{1+\epsilon}}]$,
$p_k=[ k^{-2} 2^{k^{1+\epsilon}}]$, $N_k=(m_k+p_k)k^4$.
For $\bk k\in \Bbb N^d$ and an integer $q$,
let $ \bk m_{\bk k}= (m_{k_1}, \cdots, m_{k_d})$,
$ \bk p_{\bk k}= (p_{k_1}, \cdots, p_{k_d})$,
$ \bk N_{\bk k}= (N_{k_1}, \cdots, N_{k_d})$
and $\bk k^q= (k_1^q, \cdots, k_d^q)$.
We first show the following equality
\begin{equation}\label{eq2.22}
\sum_{\bk k} \pr\big( \overline{S}_{\bk N_{\bk k}}\ge (1-7\epsilon)
\sigma (2 d |\bk N_{\bk k}| \log \log  |\bk N_{\bk k}|)^{1/2}\big)
=\infty.
\end{equation}
Set $I_{\bk i}=I_{\bk i}(\bk k) =\{ (\bk i-\bk 1)(\bk m_{\bk k}+\bk p_{\bk k})+\bk 1
\le \bk m \le (\bk i-\bk 1)(\bk m_{\bk k}+\bk p_{\bk k})+\bk m_{\bk k}\}, $
$A_{\bk k}= \bigcup_{\bk 1\le \bk i\le \bk k^4 }I_{\bk i}$ and
$B_{\bk k}= \{ \bk m: \bk m \le \bk N_{\bk k}\}\backslash A_{\bk k}$.
Then $Card A_{\bk k}=|\bk k|^4 |\bk m_{\bk k}| \sim |\bk N_{\bk k}|$
and $Card B_{\bk k}=|\bk N_{\bk k}|-Card A_{\bk k}= o(|\bk N_{\bk k}|)$.
Let
$$ v_{\bk i,1}=\sum_{\bk j\in I_{\bk i} } \overline{X}_{\bk j} $$
and
$$ S_{\bk k,1}= \sum_{\bk 1\le \bk i\le \bk k^4 } v_{\bk i,1}
=\sum_{\bk j\in A_{\bk k} } \overline{X}_{\bk j}, \quad S_{\bk
k,2}=\sum_{\bk j\in B_{\bk k} } \overline{X}_{\bk j}. $$
Clearly,
$$ \overline{S}_{ \bk N_{\bk k} }= S_{\bk k,1 }+S_{\bk k,2 }$$
and
$$ \sum_{\bk j\in B_{\bk k} } \ep \overline{X}_{\bk j}^2 / |\bk N_{\bk k}|
 \le \ep X_{\bk 1}^2 Card B_{\bk k}/|\bk N_{\bk k}| \to 0,
\quad  \bk k \to \infty. $$
From (\ref{eq2.1}), it follows that
for any $\epsilon>0$
\begin{eqnarray*}
&& \sum_{\bk k} \pr\big( |S_{\bk k,2}| \ge \epsilon
 \sigma (2 d |\bk N_{\bk k}| \log \log |\bk N_{\bk k}| )^{1/2} \big) \\
&\le& 2 \sum_{\bk k} \exp\Big\{ -\frac{ \epsilon^2 \sigma^2 2 d
|\bk N_{\bk k} | \log \log |\bk N_{\bk k} |}{
   2 \{b \epsilon \sigma (2 d |\bk N_{\bk k}| \log \log |\bk N_{\bk k}|)^{1/2}
    +\sum_{\bk j\in B_{\bk k} } \ep \overline{X}_{\bk j}^2 \} } \Big\} \\
&\le&C\sum_{\bk k} \exp\big\{ -3d \log \log |\bk N_{\bk k}| \big\}
<\infty.
\end{eqnarray*}
Thus in order to prove (\ref{eq2.22}) is enough to show that
\begin{equation}\label{eq2.23}
\sum_{\bk k} \pr\big( S_{\bk k,1} \ge
(1-6\epsilon) \sigma (2 d |\bk N_{\bk k}| \log \log |\bk N_{\bk
k}|)^{1/2} \big) =\infty.
\end{equation}
Let $B_{\bk k}^2 =\sum_{ \bk 1\le \bk i\le \bk k^4 } \ep v_{\bk i,1}^2$.
Since
$$\sum_{ \bk 1\le \bk i\le \bk k^4 } \ep\big(\sum_{\bk j\in I_{\bk i} } X_{\bk j}\big)^2
=|\bk k|^4\ep S_{\bk m_{\bk k}}^2\sim |\bk k|^4 |\bk m_{\bk k}| \sigma^2
\sim |\bk N_{\bk k}| \sigma^2$$
and
$$\sum_{ \bk 1\le \bk i\le \bk k^4 } \ep\big(v_{\bk i,1}-\sum_{\bk j\in I_{\bk i} } X_{\bk j}\big)^2
/ |\bk N_{\bk k}|\le C\ep X^2I\{|X|\ge b\}\to 0 \text{ as } b\to \infty, $$
for $b$ and $\bk k$ large enough
$$ B_{\bk k}^2 \ge (1-\epsilon)^2 |\bk N_{\bk k}| \sigma^2. $$
From Lemma \ref{lem2.3}, it
follows that
\begin{eqnarray*}
&&\pr\big(  S_{\bk k,1} \ge (1-6\epsilon) \sigma
   (2d |\bk N_{\bk k}| \log \log |\bk N_{\bk k}| )^{1/2} \big) \\
&\ge &\pr\big(  S_{\bk k,1} \ge (1-5\epsilon)
   (2d B_{\bk k}^2|\bk N_{\bk k}| \log \log |\bk N_{\bk k}| )^{1/2} \big) \\
&\ge& \big\{ 1-\Phi\big(1+(1-5\epsilon)(2 d \log \log |\bk N_{\bk
k}| )^{1/2} \big) \big\} -J_{\bk k,1}-J_{\bk k,2}
\end{eqnarray*}
where
\begin{eqnarray*}
&& J_{\bk k,1}=O(1) |\bk N_{\bk k}|^{-1} \sum_{\bk 1\le i \ne j\le
\bk k^4}
   |\ep (v_{\bk i,1} v_{\bk j,1})|, \\
&& J_{\bk k,2}=O(1) |\bk N_{\bk k}|^{-3/2} \sum_{\bk 1\le \bk i
\le \bk k^4}\ep |v_{\bk i,1}|^3.
\end{eqnarray*}
Obviously,
\begin{eqnarray*}
&& \sum_{\bk k}
\big\{ 1-\Phi\big(1+(1-5\epsilon)(2 d \log \log |\bk N_{\bk k}|)^{1/2} \big) \big\} \\
&\ge& C \sum_{\bk k} (\log |\bk N_{\bk k}|)^{-(1-4 \epsilon)d} \ge
C \sum_{\bk k} \|\bk k\|^{-(1-2\epsilon)d} =\infty,
\end{eqnarray*}
and by Lemma \ref{lem2.1},
$$ \sum_{\bk k} J_{\bk k,2}
\le C_p  \sum_{\bk k} |\bk N_{\bk k}|^{-3/2} |\bk k|^4
 ( |\bk m_{\bk k}|^{3/2} (2b)^3 +|\bk m_{\bk k}| (2b)^3)
\le C \sum_{\bk k} |\bk k|^{-2}<\infty. $$
Also,
\begin{eqnarray*}
J_{\bk k,1}&=& O(1) |\bk N_{\bk k}|^{-1}
  \sum_{\bk 1\le \bk i\ne \bk j \le \bk N_{\bk k}, \bk j-\bk i\ge \bk p_{\bk k} }
|\Cov( \overline{X}_{\bk i},  \overline{X}_{\bk j} )| \\
&=& O(1) |\bk N_{\bk k}|^{-1}
  \sum_{\bk 1\le \bk i\ne \bk j \le \bk N_{\bk k}, \bk j-\bk i\ge \bk p_{\bk k} }
|\Cov( X_{\bk i},  X_{\bk j} )| \\
&=& O(1) \sum_{\bk p_{\bk k} \le \bk j \le \bk N_{\bk k}
}|\Cov(X_{\bk 1}, X_{\bk j+\bk 1})| = O(1) \sum_{\bk N_{\bk
k-\bk1}< \bk j \le \bk N_{\bk k} }|\Cov(X_{\bk 1}, X_{\bk j+\bk
1})|,
\end{eqnarray*}
since $ p_{k_i} =o(N_{k_i})$ for  $i=1,\cdots,d$. It follows that
$$ \sum_{\bk k}J_{\bk k,1}=O(1) \sum_{\bk j} |\Cov(X_{\bk 1}, X_{\bk j+\bk 1})|<\infty. $$
Hence (\ref{eq2.23}) is proved.

Now, let $C_{\bk k} =\{ \bk m: \bk N_{\bk k-\bk 1}< \bk m\le \bk N_{\bk k}\}$,
$D_{\bk k} =\{ \bk m: \bk m\le \bk N_{\bk k}\}\backslash C_{\bk k}$ and
$$ U_{\bk k} =\sum_{ \bk j\in C_{\bk k} } \overline{X}_{\bk j} ,
\quad V_{\bk k} =\sum_{ \bk j\in D_{\bk k} } \overline{X}_{\bk j} . $$ Then
$\overline{S}_{\bk N_{\bk k}} = U_{\bk k}+V_{\bk k}$, $Card C_{\bk k}\sim
|\bk N_{\bk k}|$ and $ Card D_{\bk k} =o(|\bk N_{\bk k}|) $.
From
(\ref{eq2.1}), it follows that for any $\delta>0$,
\begin{equation}\label{eq2.24}
\sum_{\bk k} \pr\big( |V_{\bk k}| \ge \delta \sigma
 (2 d |\bk N_{\bk k}| \log \log |\bk N_{\bk k}|)^{1/2} \big)<\infty.
 \end{equation}
By the Borel-Cantelli lemma, we conclude that
\begin{equation}\label{eq2.25} \limsup_{\bk k\to \infty} \frac {|V_{\bk k}|}{(2
   d |\bk N_{\bk k}| \log \log |\bk N_{\bk k}|)^{1/2} }= 0\quad a.s.
\end{equation}
 From (\ref{eq2.22}) and (\ref{eq2.24}), it follows that
\begin{equation}\label{eq2.26}
\sum_{\bk k} \pr\big( U_{\bk k} \ge (1-8\epsilon) \sigma
 (2 d |\bk N_{\bk k}| \log \log |\bk N_{\bk k}|)^{1/2} \big)=\infty.
 \end{equation}
Note that $\{ U_{\bk k}; \bk k \in \Bbb N^d\}$ is a negatively associated filed. It follows that
for any $x, y $ and $ \bk i\ne \bk j$,
$$ \pr (U_{\bk i}\ge x, U_{\bk j}\ge y)\le \pr (U_{\bk i}\ge x)\pr (U_{\bk j}\ge y). $$
Hence, by the generalized Borel-Cantelli lemma, (\ref{eq2.26})
yields
\begin{equation}\label{eq2.27}
\limsup_{\bk k\to \infty} \frac { U_{\bk k} }{(2
   d |\bk N_{\bk k}| \log \log |\bk N_{\bk k}|)^{1/2} } \ge (1-8\epsilon) \sigma
\quad a.s.
\end{equation}
From (\ref{eq2.25}) and (\ref{eq2.27}), it follows that
(\ref{eq2.21}) holds.

\bigskip
{\noindent \bf Proof of Theorem \ref{th1}} is similar to that of (\ref{upbound})
in which we choose $m=1$ and $Y_{\bk k}=\overline{X}_{\bk k}$ instead.

\bigskip
{\noindent \bf Proof of Theorem \ref{th3}:} From  (\ref{eq1.6}),
it follows that there exists a constant $0<C<\infty$ such that
$$\pr\left( \limsup_{\bk n\to \infty} \frac{ |S_{\bk n}|}{ (2 d |\bk n| \log\log |\bk n|)^{1/2} }<
C \right)>0. $$
Then
\begin{equation}\label{eq2.28}
\pr\left( \limsup_{\bk n\to \infty} \frac{  |X_{\bk n}|}{
  (2 d |\bk n| \log\log |\bk n|)^{1/2} }< 2C \right)>0.
  \end{equation}
Let
\begin{eqnarray*}
&&A_{\bk n}=  \{|X_{\bk n}|\ge 2C  (2 d |\bk n| \log\log |\bk n|)^{1/2} \}, \\
&&A_{\bk n}^+ =  \{ X_{\bk n}^+ \ge 2C  (2 d |\bk n| \log\log |\bk n|)^{1/2} \}, \\
&& A_{\bk n}^- =  \{ X_{\bk n}^- \ge 2C  (2 d |\bk n| \log\log
|\bk n|)^{1/2} \}.
\end{eqnarray*}
Note that $\big\{ I\{A_{\bk n}^+\}; \bk n\in \Bbb N^d\big\}$ and
$\big\{ I\{A_{\bk n}^-\}; \bk n\in \Bbb N^d\}$
are both negatively associated fields. It follows that for any $\bk m$ and $\bk n$,
\begin{eqnarray*}
&& \Var \Big\{ \sum_{\bk m\le \bk k \le \bk m+\bk n}I\{A_{\bk
n}^+\}\Big\}
\le \sum_{\bk m\le \bk k \le \bk m+\bk n}\pr (A_{\bk n}^+), \\
&& \Var \Big\{ \sum_{\bk m\le \bk k \le \bk m+\bk n}I\{A_{\bk
n}^-\}\Big\} \le \sum_{\bk m\le \bk k \le \bk m+\bk n}\pr (A_{\bk
n}^-).
\end{eqnarray*}
It follows that
$$ \Var \Big\{ \sum_{\bk m\le \bk k \le \bk m+\bk n}I\{A_{\bk k}\}\Big\}
\le 2 \sum_{\bk m\le \bk k \le \bk m+\bk n}\pr (A_{\bk k}). $$
Hence from  Lemma A.6 of Zhang and Wen (2001b), it follows that
for any $\bk m$
\begin{equation}\label{eq2.29}
\big( 1- \pr( \bigcup_{ \bk k\ge \bk m} A_{\bk k} ) \big)^2
\sum_{\bk k\ge \bk m} \pr( A_{\bk k})
\le 2 \pr( \bigcup_{ \bk k\ge \bk m} A_{\bk k} ).
\end{equation}
Since (\ref{eq2.28}) implies $ \pr\big( A_{\bk n}, i.o.\big) <1$,
we conclude that for some $\bk m$,
$$ \pr( \bigcup_{ \bk k\ge \bk m} A_{\bk k} ):=\beta<1. $$
Then from (\ref{eq2.29}), it follows that
$$\sum_{\bk k\ge \bk m} \pr( A_{\bk k})
\le \frac{ 2\beta}{(1-\beta^2)}<\infty. $$
So,
$$ \sum_{\bk k\ge \bk m}
\pr\big( |X_{\bk 1}|\ge 2C  (2 d |\bk k| \log\log |\bk k|)^{1/2}\big)<\infty, $$
which implies
\begin{equation}\label{eq2.30}
\ep X_{\bk 1}^2 \log^{d-1}( |X_{\bk 1}|) / \log\log (|X_{\bk
1}|)<\infty.
\end{equation}
 Finally, from (\ref{eq2.30}) and
the law of the large numbers (c.f. Zhang and Wang 1999), it
follows that
$$ \lim_{\bk n} \frac{ S_{\bk n}- |\bk n| \ep X_{\bk 1}}{ |\bk n| } \to 0 \quad a.s. $$
which together with (\ref{eq1.6}) yields $\ep X_{\bk 1}=0$. And
then (\ref{eq1.4}) holds.

\bigskip

%\bigskip{\bf Acknowledgements:} The author is very grateful to the referee for his or her careful reading of the manuscript.

%\newpage

\bigskip
\bigskip
\bigskip

\noindent{\sc Li-Xin Zhang\\
Department of Mathematics  \\
Zhejiang University, Yuquan Campus\\
Hangzhou 310027, Zhejiang\\
P.R. China\\
Email:  stazlx@zju.edu.hz.cn }

\end{document}